\pdfoutput=1

\documentclass{amsart}

\usepackage{my_style}
\usepackage{paralist} 

\usepackage{tikz}
\usepackage{tikz-3dplot}
\usepackage{tikz-cd}
\usetikzlibrary{calc}

\tikzset{%
	my arrow/.style = {thick,->, color=blue},
	my triangle/.style = {color=blue,fill opacity=0.7,fill=green!80!blue, thick},
	big tetrahedron/.style = {thick, opacity=0.4},
	small tetrahedron/.style = {color=blue,fill opacity=0.4,fill=green!80!blue, thick, opacity=0.4},
	my axes/.style = {ultra thin}
}

\title[Monomial Cremona transformations of the plane]{Classification of the monomial Cremona transformations of the plane}

\author{Corey Harris}
\address{Department of Mathematics, Florida State University, Tallahassee FL, 32306, USA}
\email{charris@math.fsu.edu}
\urladdr{\url{http://www.coreyharris.name}}

\newcommand{\ch}{\operatorname{ch}}
\newcommand{\cone}{\operatorname{cone}}
\newcommand{\Vol}{\operatorname{Vol}}
\newcommand{\lr}[2]{$\langle #1 ; \; #2 \rangle$}

\begin{document}

\begin{abstract}
We classify all monomial planar Cremona maps by multidegree using recent methods developed by Aluffi. Following the main result, we prove several more properties of the set of these maps, and also extend the results to the more general `r.c. monomial' maps.
\end{abstract}

\maketitle

\section{Introduction}  \label{sec:intro}

Let $(x:y:z)$ be coordinates on $\bbP^2$.
Let $\varphi: \bbP^2 \dashrightarrow \bbP^2$ be a monomial map 
\[ \varphi: (x:y:z) \mapsto (x^{a_{11}} y^{a_{12}} z^{a_{13}}: x^{a_{21}} y^{a_{22}} z^{a_{23}}: x^{a_{31}} y^{a_{32}} z^{a_{33}}).\]
We can carry the information of this map in its \emph{exponent matrix} 
$M_\varphi = (a_{ij})$.
We will be interested in birational monomial maps on $\bbP^2$.  For such maps, 
the total degree of each monomial is constant,
i.e., there is a $\delta$ such that the row sum is $\delta$ for each row of $M_\varphi$.

If the monomials $x^{a_{i1}} y^{a_{i2}} z^{a_{i3}}$ share no common factors, we say that 
$\varphi$ is written in \emph{reduced form}.  Throughout the paper, we'll assume that any 
monomial map $\bbP^2 \dashrightarrow \bbP^2$ comes in reduced form. Of course, if $\varphi$ 
is in reduced form, then $M_\varphi$ must have a 0 in each column.  Thus, if $M_\varphi$ is 
the exponent matrix of a monomial Cremona transformation on $\bbP^2$, then up to swapping of 
rows and columns it has one of the forms
\[
\left( \begin{matrix}
0 & 0 & \cdot \\
\cdot & \cdot & 0 \\
\cdot & \cdot & \cdot
\end{matrix} \right)
\text{ or }
\left( \begin{matrix}
0 & \cdot & \cdot \\
\cdot & 0 & \cdot \\
\cdot & \cdot & 0
\end{matrix} \right).
\]
A rational map on the projective plane has associated to it a tuple of numbers 
$(\gamma_0,\gamma_1,\gamma_2)$ called the multidegree (see 
\autoref{sec:multidegrees_via_polyhedra}).  A rational map is a Cremona map 
if and only if its multidegree is $(1,d,1)$, and if the map is monomial then 
$d = \delta$ as above.

\autoref{thm:mainThm} below gives the complete list of monomial Cremona 
transformations of the plane.  Notice that (\bf II\rm) could actually have been 
included in (\bf III\rm) if we just allowed $c=0$.  We list it separately to 
emphasize that it is special in the sense that the exponent matrix has more than 
3 zeroes.

\begin{thm} \label{thm:mainThm}
Let $\varphi: \bbP^2 \rightarrow \bbP^2$ be a monomial rational map with exponent 
matrix $M_\varphi$.  Then $\varphi$ is a Cremona transformation (with multidegree 
$(1,\delta,1)$) if and only if 
\begin{enumerate}[\rm(\bf I\rm)]
	\item $\varphi$ is the standard involution 
		$\varphi(x:y:z) \mapsto (xy:xz:yz)$, or
	\item $M_\varphi$ is of the form
		\[ M_\varphi = \left( \begin{matrix}
			0 & 0 & \delta\\
			1 & \delta-1 & 0\\
			0 & 1 & \delta -1
		\end{matrix} \right), \]
		or
	\item $M_\varphi$ is of the form
		\[ M_\varphi = \left(
		\begin{matrix}
			0 & 0 & \delta \\
			a & b & 0 \\
			c & d & e
		\end{matrix}
		\right) \]
		and the following equations are satisfied:
		\begin{enumerate}[\rm(i)]
			\item $a+b = \delta$,\; $c+d+e=\delta$,
			\item $ad-bc=1$.
		\end{enumerate}
\end{enumerate}
\end{thm}

\begin{remark}
Throughout the paper we will allow ourselves to swap two rows or columns of 
$M_\varphi$ whenever convenient, as we've done here.
\end{remark}

A useful fact about $M_\varphi$ is the
\begin{lemma} \label{lem:determinant}
$|\det(M_\varphi)| = \delta$.
\begin{proof}
See \cite[Prop. 3.1]{GSP03}, \cite[Prop. 1]{DL14}, \cite[Sec. 2]{Joh14}.
\end{proof}
\end{lemma}


\section{Some convex geometry} \label{sec:convex_geom}

The rows of $M_\varphi$ give coordinates for points in $\bbR^3$.
\begin{defn}
	Let $S = \setrow sk$ be a set of points in $\bbR^n$.  We denote the 
	\emph{convex hull} of $S$ by \[\ch(S) := 
	\{ \lambda_1 s_1 + \dots + \lambda_k s_k \;|\; \lambda_i \geq 0, \lambda_1 + \dots + \lambda_k = 1 \}\]
	and the \emph{conical hull} of $S$ by \[\cone(S) := 
	\{ \lambda_1 s_1 + \dots + \lambda_n s_k  \;|\; \lambda_i \geq 0 \}.\]
	In general, if $R$ is an arbitrary subset of $\bbR^n$, we define its convex hull to be
	\[ ch(R) := \bigcup_{a,b \in R} \; ch(\{a,b\}) \]
	Here we are concerned with \emph{convex polyhedra}, subsets 
	$P \subset \bbR^n$ which can be written \[ P = \ch(S) + \cone(S') := 
	\{ a + b \;|\; a \in ch(S), \; b \in cone(S') \}. \]
	For a polyhedron of dimension $d$, the faces of dimension 0,1 and 
	$d-1$ are called \emph{vertices}, \emph{edges} and
	\emph{facets}, respectively.
\end{defn}

Let $\mathcal{B} := \{e_1, e_2, e_3\}$ be the standard basis in $\bbR^3$.  
If 
\[ M_\varphi = \left( \begin{matrix}
A \\
B \\
C
\end{matrix} \right) \]
is the exponent matrix of a monomial rational map $\varphi$,
then the \emph{Newton polyhedron} of $\varphi$ is 
$N = \ch(\{A,B,C\}) + \cone(\mathcal{B})$.  This is a 3-dimensional, unbounded 
polyhedron with exactly one finite facet: $\ch(\{A,B,C\})$.

We will adopt the following notation for a convex polyhedron.  If 
$S = \setrow vj$ and $S' = \setrow wk$ then we write
\lr{\row vj}{\row wk} to mean $\ch(S) + \cone(S')$.
We say that the \emph{order} of the polyhedron is the minimal number $|S|+|S'|$.
Thus, the Newton polyhedron $N$ is \lr{A, B, C}{e_1, e_2, e_3}.

\subsection{Triangulations} \label{sec:triangulations}
We will need to triangulate $N$.  In this paper, a triangulation of a 3-dimensional convex polyhedron $N$ is a set 
$\mathcal{T}$ of 3-dimensional subsets of $N$ such that 
$\sqcup_{P \in \mathcal{T}} P = N$ and such that the order of $P \in \mathcal{T}$ 
is 4 for all $P$ (this makes $P$ ``simplicial'').

A useful algorithm for doing this comes from \cite{Cla85}.  We choose a 
distinguished vertex $p \in \{A,B,C\}$ and construct a triangulation 
$\mathcal{T}$ starting with $\langle p; e_1, e_2, e_3 \rangle$, that is, 
the positive orthant translated to $p$. Now let $F_p$ denote the facets of 
$N$ which do not include the vertex $p$.  If $f \in F_p$ is a facet with 
order 3, then include $\ch(\{p\} \cup f)$ in $\mathcal{T}$.  Otherwise, we 
should construct a triangulation $\mathcal{T}_f$ of $f$ using the analogous 
procedure and for each $t \in \mathcal{T}_f$, include $\ch(\{p\} \cup f)$ 
in $\mathcal{T}$.

\tdplotsetmaincoords{65}{110}

\begin{figure}[hb]
	\centering
	\begin{tikzpicture}[scale=2,tdplot_main_coords]
	    \coordinate (O) at (0,0,0);
	    \coordinate (e1) at (1,0,0);
	    \coordinate (e2) at (0,1,0);
	    \coordinate (e3) at (0,0,1);
	    \coordinate (A1) at (0,0,1);
		\coordinate (A2) at (0,1,0);
		\coordinate (A3) at (0.65,0.35,0);
	    \draw [my axes] (O) -- ($ (e1) + (e1) $) node[anchor=north east]{$x$};
	    \draw [my axes] (O) -- ($ (e2) + (e2) $) node[anchor=west]{$y$};
		\draw [my axes] (O) -- ($ (e3) + (e3) $) node[anchor=south]{$z$};
		\draw[my triangle] (A1) -- (A2) -- (A3) -- cycle;
		\draw[my arrow] (A1) -- ($ (A1) + (0,0,0.8) $);
		\draw[my arrow] (A2) -- ($ (A2) + (0,0.8,0) $);
		\draw[my arrow] (A3) -- ($ (A3) + (e1) $);
		\draw[my arrow] (A1) -- ($ (A1) + (e1) $);
		\node [above left] at (A1) {A};
		\node [above right] at (A2) {B};
		\node [below right] at (A3) {C};
	\end{tikzpicture}
	\caption{Schematic drawing of a Newton polyhedron}
	\label{fig:triangulationPolyhedron}
\end{figure}
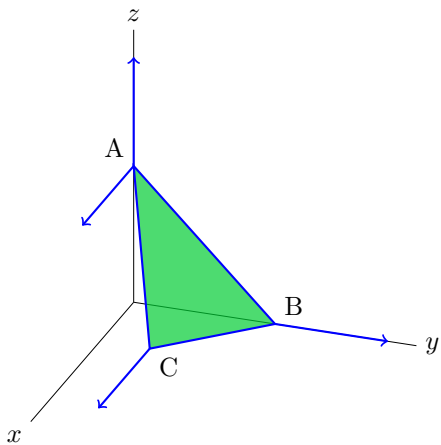

\begin{ex} 
	In \autoref{fig:triangulationPolyhedron}, where $A,B,C$ are given by the matrix 
	\[ M_\varphi = \left( \begin{matrix}
	0 & 0 & 3 \\
	0 & 3 & 0 \\
	2 & 1 & 0
	\end{matrix} \right), \]
	
	we can begin a triangulation $\mathcal{T}$ of $N$ with 
	$\langle A; e_1, e_2, e_3 \rangle$. We then consider the facets which do not include $A$.  The only such face is $\langle B,C; e_1,e_2 \rangle$. 
	Since this face has order $4 > 3$, we must triangulate it.  We choose a 
	vertex, say $B$, and take the cone $\langle B; e_1,e_2 \rangle$ at $B$. 
	The remaining facet (of $\langle B,C; e_1,e_2 \rangle$) is 
	$\langle C; e_2 \rangle$.  Thus, the two polyhedra 
	$\langle B; e_1,e_2 \rangle$ and $\langle B,C;e_2 \rangle$ form a 
	triangulation of $\langle B,C; e_1,e_2 \rangle$.  Taking the pyramid of 
	these polyhedra at $A$ completes the triangulation of $N$:
	\[ \mathcal{T} = \{ 
	\langle A;e_1,e_2,e_3 \rangle, 
	\langle A,B; e_1,e_2 \rangle, 
	\langle A,B,C; e_2 \rangle \}. 
	\]
\end{ex}

\section{Computing multidegrees via polyhedra} \label{sec:multidegrees_via_polyhedra}


Let $\Gamma(\varphi)$ be the closure of the graph of $\varphi$ in $\bbP^2 \times \bbP^2$.
\begin{figure}[t]
	\centering
	\begin{tikzcd}[ampersand replacement=\&]
		\& \Gamma(\varphi) \subset \bbP^2 \times \bbP^2 \arrow{dl} \arrow{dr} \& \\
		\bbP^2 \arrow{rr}{\varphi} \& \& \bbP^2
	\end{tikzcd}
\end{figure}
Then the class $[\Gamma(\varphi)] \in A(\bbP^2 \times \bbP^2)$ is 
\[ 	[\Gamma(f)] = \gamma_0 h_2^2 + \gamma_1 h_1 h_2 + \gamma_2 h_1^2 \]
where $h_1,h_2$ are the pullbacks of the hyperplane class $c_1(\calO(1))$ from $\bbP^2$ via 
the projections to the first and second factors, respectively.

\begin{defn}
The \emph{multidegree} of the rational map $\varphi: \bbP^2 \dashrightarrow \bbP^2$ is the tuple $(\gamma_0,\gamma_1,\gamma_2)$.
\end{defn}

The coefficients $\gamma_i$ satisfy
\begin{enumerate}[(\bf i\rm)]
\item $\gamma_0 = \gamma_0(\varphi) := \# (\varphi^{-1}(p))$ for a general point $p \in \bbP^2$,
\item $\gamma_1 = \gamma_1(\varphi)$ is the degree of (the closure of) 
$\varphi^{-1}(\bbP^1) \subset \bbP^2$ for a general hyperplane $H \subset \bbP^2$,
\item $\gamma_2 = \gamma_2(\varphi)$ is the degree of the field extension 
$K(\varphi(\bbP^2)) \subset K(\bbP^2)$.
\end{enumerate}
If $\varphi$ is a monomial Cremona transformation, then $\gamma_0 = 1$ because 
$\varphi$ is generically one-to-one, and $\gamma_2=1$ because $\varphi$ is 
birational.  The degree of the image under $\varphi$ of a hyperplane in $\bbP^2$ 
will be $\# (H \cdot \varphi^{-1} (H') ) = \delta$, the total degree of the 
monomials defining $\varphi$.  The conclusion is that a monomial rational map 
$\varphi$ on $\bbP^2$ is a Cremona transformation if and only if the 
multidegree of $\varphi$ is $(1,\delta,1)$.

\subsection{}
Let $\mathcal{T}$ be a triangulation of $N$, the Newton polyhedron of $\varphi$.
If $\langle S;V \rangle \in \mathcal{T}$, we define $\pi_V$ to be 
the projection $\bbR^3 \rightarrow V^\perp$ where $V^\perp$ denotes the subspace 
spanned by $\mathcal{B} \backslash V$.  We then say (in the language of 
\cite[Sec. 3.1]{Alua}) that the \emph{volume} $\Vol(\langle S;V \rangle)$ 
is the normalized volume of $\pi_V(\langle S; V \rangle)$.

Our method for computing multidegrees is \cite[Thm. 1.4]{Alua},
which we restate here using our notation.  Letting 
$T(i) = \{ \langle S;V \rangle \in \mathcal{T} \;|\; \#(S) = i+1 \}$ be 
the polyhedra $\langle S;V \rangle$ whose projection under $\pi_V$ are of 
dimension $i$, we have
\begin{thm}[Aluffi]
\label{thm:AluffiMain}
$\gamma_i = \sum_{P \in \mathcal{T}(i)} \Vol(P)$.
\end{thm}
Via this theorem, we have the following strategy.  Find a triangulation of the 
Newton polyhedron $N$ associated to an exponent matrix $M_\varphi$.  If the 
volumes in the triangulation give $(\gamma_0,\gamma_1,\gamma_2) = (1, \delta, 1)$,
then $\varphi$ is a Cremona transformation, otherwise not.

\section{Main analysis} \label{sec:main_analysis}

In this section, we use the technology outlined in the previous sections to 
analyze the types of monomial Cremona transformations that exist, 
resulting in the proof of \autoref{thm:mainThm}.  The 
specific problem is: given a $3 \times 3$ integer matrix, when is it the (reduced)
exponent matrix of a monomial planar Cremona transformation?

\subsection{Main case} \label{subsec:main_case}
We begin with the case where the exponent matrix 
$ M_\varphi = \left(\begin{smallmatrix} A \\ B \\ C \end{smallmatrix}\right) $
has the form
\[ \left(\begin{matrix}
0 & 0 & \cdot \\
\cdot & \cdot & 0 \\
\cdot & \cdot & \cdot
\end{matrix} \right)
= \left( \begin{matrix}
0 & 0 & \delta \\
a & b & 0 \\
c & d & e
\end{matrix} \right). \]
Note, we assume that $c,d$ are both non-zero, since otherwise we 
would be in the other case (cf. \autoref{sec:intro}).

If $M_\varphi$ is going to be an exponent matrix for a monomial Cremona transformation then 
we must have
\begin{equation}
\label{eqn:main_case_row_sum}
a+b = \delta = c+d+e
\end{equation}
and by \autoref{lem:determinant} we must also have $|\delta(ad-bc)| = \delta$, 
so $ad-bc = \pm 1$.  
By swapping columns 1 and 2 if necessary, we may assume
\begin{equation}
\label{eqn:main_case_det}
ad-bc = 1
\end{equation}
which uniquely determines the edges of the Newton polyhedron $N$ associated to $M_\varphi$ 
(consider the projection of the points $B,C$ to the $xy$-plane), 
an example of which is shown in \autoref{fig:mainCase}.

\begin{figure}[bth]
	\centering
	\begin{tikzpicture}[scale=.45,tdplot_main_coords]
	    \coordinate (O) at (0,0,0);
	    \coordinate (e1) at (5,0,0);
	    \coordinate (e2) at (0,5,0);
	    \coordinate (e3) at (0,0,5);
	    \coordinate (A1) at (0,0,5);
		\coordinate (A2) at (2.1,2.9,0);
		\coordinate (A3) at (2,1,2);
	    \draw [my axes] (O) -- ($ (e1) + (e1) $) node[anchor=north east]{$x$};
	    \draw [my axes] (O) -- (0,8,0) node[anchor=west]{$y$};
		\draw [my axes] (O) -- ($ (e3) + (e3) $) node[anchor=south]{$z$};
		\draw[my triangle] (A1) -- (A2) -- (A3) -- cycle;
		\draw[my arrow] (A1) -- ($ (A1) + (0,0,4) $);
		\draw[my arrow] (A1) -- ($ (A1) + (0,4,0) $);
		\draw[my arrow] (A1) -- ($ (A1) + (e1) $);
		\draw[my arrow] (A2) -- ($ (A2) + (0,4,0) $);
		\draw[my arrow] (A2) -- ($ (A2) + (e1) $);
		\draw[my arrow] (A3) -- ($ (A3) + (e1) $);
		\node [above left] at (A1) {A};
		\node [above right] at (A2) {B};
		\node [left] at (A3) {C};
	\end{tikzpicture}
	\caption{}
	\label{fig:mainCase}
\end{figure}

Begin constructing the triangulation $\mathcal{T}$ with $\langle B; e_1,e_2,e_3 \rangle$.  The 
facets of $N$ which do not contain $B$ are
\[ 
\langle A; \; e_1,e_3 \rangle, \;
\langle A; \; e_2,e_3 \rangle, \;
\langle A,C; \; e_1 \rangle, 
\]
so we have
\[ 
\mathcal{T} = \{ 
\langle B; \; e_1,e_2,e_3 \rangle, \;
\langle A,B; \; e_1,e_3 \rangle, \;
\langle A,B; \; e_2,e_3 \rangle, \;
\langle A,B,C; \; e_1 \rangle
\}.
\]
The polyhedron $\langle B; e_1,e_2,e_3 \rangle$ projects to the point $B$, which 
has $Vol(B) = 1 = \gamma_0$.  
The projection of $\langle A,B; e_2,e_3 \rangle$ onto the $x$-axis is the 
interval $(0,a)$, which has $\Vol((0,a)) = a$.
The projection of $\langle A,B; e_1,e_3 \rangle$ onto the $y$-axis is the 
interval $(0,b)$, which has $\Vol((0,b)) = b$ and as we would hope, we get 
$\gamma_1 = a + b = \delta$.

Finally, we consider $\langle A,B,C; e_1 \rangle$ which should be projected 
onto the $yz$-plane. The projection is a triangle with vertices
$(0,\delta),\; (b,0),\; (d,e)$, and its volume is
\[
	\Vol(\langle A,B,C; e_1) = \big| (2!) (\frac{1}{2!})
	\left| \begin{matrix}
	0 & \delta & 1 \\
	b & 0 & 1 \\
	d & e & 1
	\end{matrix} \right| \big|
	= |\delta d + be - \delta b|
\]
So since we should have $\gamma_2 = 1$ we have the necessary condition
\begin{equation}
	\label{eqn:main_case_gamma1}
	|\delta d + be - \delta b| = 1.
\end{equation}
Notice though, that if we use $a = \delta -b$ and $c = \delta - d - e$ in the determinant equation $ad-bc=1$, we get
\[ (\delta - b)d - b(\delta - d -e) = \delta d +be - \delta b = 1. \]

This shows that \autoref{eqn:main_case_gamma1} is equivalent to the requirement 
that $ad-bc=1$ as long as we require also that \autoref{eqn:main_case_row_sum} 
be satisfied, completing the proof of the
\begin{lemma}
\label{lemma:mainCase}
If $c >0, d>0$, then the exponent 
matrix 
\[ M_\varphi = \left(
\begin{matrix}
0 & 0 & \delta \\
a & b & 0 \\
c & d & e
\end{matrix}
\right) \]
defines a monomial Cremona transformation with multidegree $(1,\delta,1)$ if and
only if the following equations are satisfied:
\begin{enumerate}[\rm(\bf i\rm)]
\item $a+b = \delta$,\; $c+d+e=\delta$,
\item $ad-bc=1$.
\end{enumerate}
\end{lemma}


For a quick payoff, if we set $a = \delta-1$, then we find a uniquely
determined exponent matrix:
\[ M_\varphi = \left( 
\begin{matrix}
0 & 0 & \delta \\
\delta-1 & 1 & 0 \\
\delta-2 & 1 & 1
\end{matrix}
\right) \]
and one can check that this yields a monomial Cremona transformation with 
inverse given by the exponent matrix
\[ 
	M_{\varphi^{-1}} = \left( \begin{matrix}
		1 & \delta-1 & 0 \\
		0 & 0 & \delta \\
		1 & \delta-2 & 1
	\end{matrix} \right). 
\]

\subsection{Other cases} \label{subsec:other_cases}
Now assume $M_\varphi$ has the form
\[ \left( \begin{matrix}
0 & \cdot & \cdot \\
\cdot & 0 & \cdot \\
\cdot & \cdot & 0
\end{matrix} \right) 
= \left(\begin{matrix}
0 & a & b \\
c & 0 & d \\
e & f & 0
\end{matrix} \right). \]
Here we have $a,b,c,d,e,f \in \bbZ_{\geq 0}$ with
\begin{equation}
\label{eqn:rowsums2}
	a+b = c+d = e+f = \delta.
\end{equation}

First note that \autoref{lem:determinant} in this case gives
$|ade + bcf| = \delta$.  Rewriting using \eqref{eqn:rowsums2}, we get 
\[ ade+bcf = a(\delta-c)e + (\delta-a)c(\delta-e) = \delta(\delta c + ae - ac - ce) \]
so 
\begin{equation} 
\label{eqn:determinant2}
	|\delta c + ae - ac -ce| = 1. 
\end{equation}

\renewcommand{\thefootnote}{\fnsymbol{footnote}}

\subsubsection{No additional zeroes} 
\begin{figure}[htb]
\centering
	\begin{tikzpicture}[scale=2,tdplot_main_coords]
	    \coordinate (O) at (0,0,0);
	    \coordinate (A1) at (0,1,1);
		\coordinate (A2) at (1,0,1);
		\coordinate (A3) at (1,1,0);
		\coordinate (e1) at (1,0,0);
	    \coordinate (e2) at (0,1,0);
	    \coordinate (e3) at (0,0,1);
	    \draw [my axes] (0,0,0) -- (2,0,0) node[anchor=north east]{$x$};
	    \draw [my axes] (0,0,0) -- (0,2,0) node[anchor=west]{$y$};
		\draw [my axes] (0,0,0) -- (0,0,2) node[anchor=south]{$z$};
		\draw[my triangle] (A1) -- (A2) -- (A3) -- cycle;
		\draw[my arrow] (0,1,1) -- (0,2,1);
		\draw[my arrow] (0,1,1) -- (0,1,2);
		\draw[my arrow] (1,0,1) -- (2,0,1);
		\draw[my arrow] (1,0,1) -- (1,0,2);
		\draw[my arrow] (1,1,0) -- (2,1,0);	
		\draw[my arrow] (1,1,0) -- (1,2,0);	
		\node [above right] at (A1) {A};
		\node [left] at (A2) {B};
		\node [below] at (A3) {C};
	\end{tikzpicture}
	\caption{Newton polyhedron for the standard involution}
	\label{fig:involution_poly}
\end{figure}
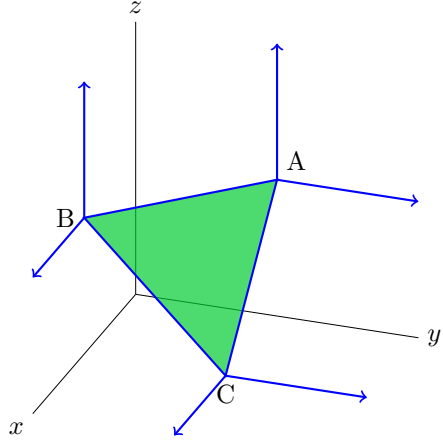

Assume $a,b,c,d,e,f$ are all nonzero.  Then $A$ is the unique minimal vertex for $e_1$\footnote{By `minimal vertex for $v$', we mean $x \cdot v$ is minimized by the row vector $x$.},
$B$ is the unique minimal vertex for $e_2$, and 
$C$ the unique minimal vertex for $e_3$.  This determines the ridge structure 
of the Newton polyhedron $N$.  An example of this structure is shown in 
\autoref{fig:involution_poly}, where 
${\varphi(x:y:z)=(yz:xz:xy)}$.

If we begin our triangulation with the distinguished point $A$, we get
\[ \mathcal{T} = \{ 
\langle A; e_1,e_2,e_3 \rangle,
\langle A, B; e_1,e_3 \rangle,
\langle A, C; e_2, e_3 \rangle,
\langle A,B,C; e_1 \rangle
\}. \]
The projection of $A$ along all three directions $e_i$ is a point with volume 
$\gamma_0 = 1$.  The projection of $\ch(A,B)$ along $e_1,e_3$ is the interval 
$(0,a)$, and the projection of $\ch(A,C)$ along $e_2,e_3$ is the interval 
$(0,e)$.  The sum of these two volumes is $\gamma_1 = a + e$.  We want 
$\gamma_1 = \delta$, so this gives a new condition: 
\begin{equation}
\label{eqn:new_equalities}
b=e \text{ and } a = f.
\end{equation}

The projection of $\ch(A,B,C)$ onto the $y z$-plane has normalized volume
\[ \big| \det \left(\begin{matrix}
a & b & 1 \\
0 & d & 1 \\
f & 0 & 1
\end{matrix} \right) \big|
= |ad + bf - df| \]
and we want $\gamma_2 = 1$ so we should require 
\begin{equation}
|ad+bf-df|=1
\end{equation}
which becomes $|ab|=1$ by \eqref{eqn:new_equalities}.  
Then, $a=1=b$ and so $e=1=f$ by \eqref{eqn:new_equalities} and $c=1=d$
by \eqref{eqn:determinant2}.
Thus, 
\[ M_\varphi = \left( 
\begin{matrix}
0 & 1 & 1 \\
1 & 0 & 1 \\
1 & 1 & 0
\end{matrix}
\right). \]

\subsubsection{$f=0$}
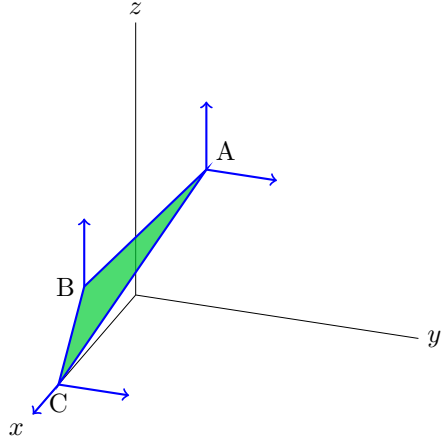
\begin{figure}[htb]
\centering
	\begin{tikzpicture}[scale=1,tdplot_main_coords]
	    \coordinate (O) at (0,0,0);
	    \coordinate (A1) at (0,1,2);
		\coordinate (A2) at (2,0,1);
		\coordinate (A3) at (3,0,0);
		\coordinate (e1) at (4,0,0);
	    \coordinate (e2) at (0,4,0);
	    \coordinate (e3) at (0,0,4);
	    \draw [my axes] (0,0,0) -- (4,0,0) node[anchor=north east]{$x$};
	    \draw [my axes] (0,0,0) -- (0,4,0) node[anchor=west]{$y$};
		\draw [my axes] (0,0,0) -- (0,0,4) node[anchor=south]{$z$};
		\draw[my triangle] (A1) -- (A2) -- (A3) -- cycle;
		\draw[my arrow] (A1) -- (0,1,3);
		\draw[my arrow] (A1) -- (0,2,2);
		\draw[my arrow] (A3) -- (3,1,0);
		\draw[my arrow] (A3) -- (4,0,0);
		\draw[my arrow] (A2) -- (2,0,2);
		\node [above right] at (A1) {A};
		\node [left] at (A2) {B};
		\node [below] at (A3) {C};
	\end{tikzpicture}
	\caption{Newton polyhedron when $f=0$}
	\label{fig:other_poly}
\end{figure}

If $f=0$, then $e = \delta$, so we have
$|ade+bcf| = |\delta a d| = \delta$ which implies $a=1, d=1$.
Then $b = \delta-1$ and $c=\delta-1$, so
\[ M_\varphi = \left( 
\begin{matrix}
0 & 1 & \delta-1 \\
\delta-1 & 0 & 1 \\
\delta & 0 & 0
\end{matrix}
\right). \]
The Newton polyhedron in this case is shown in \autoref{fig:other_poly}.

Our previous triangulation $\mathcal{T}$ is still a triangulation here.  
The facet $\langle B,C; e_1,e_3 \rangle$ can be triangulated with
\[ \{ \langle B; e_1,e_3 \rangle, \langle B,C; e_1 \rangle \} \]
which were the facets we used before, so the algorithm still generates
\[ \mathcal{T} = \{ 
\langle A; e_1,e_2,e_3 \rangle,
\langle A, B; e_1,e_3 \rangle,
\langle A, C; e_2, e_3 \rangle,
\langle A,B,C; e_1 \rangle
\}. \]
Then to get a Cremona map, we should still require $|ad+bf-df|=1$, 
and indeed this becomes
\[ |1\cdot 1 + 0 - 0| = 1. \]
This also takes care of the case where one of $a,c,d$ is zero, since we can take 
the resulting matrix to one of the proper form by row and column swaps.  For instance, if $d=0$, we get
\[\left(\begin{matrix}
0 & a & b \\
c & 0 & d \\
e & f & 0
\end{matrix} \right)
=
\left(\begin{matrix}
0 & a & b \\
c & 0 & 0 \\
e & f & 0
\end{matrix} \right)
\leadsto
\left(\begin{matrix}
0 & b & a \\
e & 0 & f \\
c & 0 & 0
\end{matrix} \right) \]
by swapping rows 2,3 and swapping columns 2,3.

\subsubsection{$e=0$}
If $e=0$, then $f=\delta$, so we have
$|ade+bcf| = |\delta bc|=\delta$, so $b=1,c=1$.
Then $a = \delta-1$ and $d=\delta-1$, so 
\[ M_\varphi = \left( 
\begin{matrix}
0 & \delta-1 & 1 \\
1 & 0 & \delta-1 \\
0& \delta & 0
\end{matrix}
\right). \]
But this is the same as in the previous case, just swap rows 1,2 and columns 1,2.

So this section concludes the analysis of exponent matrices which are not of the form of \autoref{lemma:mainCase}.  The result is thus
\begin{lemma}
\label{lemma:otherCase}
Let $\varphi: \bbP^2 \rightarrow \bbP^2$ be a monomial rational map where 
$M_\varphi$ is not in the form of \autoref{lemma:mainCase}.  Then $\varphi$ is 
a Cremona transformation (with multidegree $(1,\delta,1)$) if and only if 
\[ M_\varphi = \left( 
\begin{matrix}
0 & \delta-1 & 1 \\
1 & 0 & \delta-1 \\
0& \delta & 0
\end{matrix}
\right) 
\text{ or }\;
M_\varphi = \left( 
\begin{matrix}
0 & 1 & 1 \\
1 & 0 & 1 \\
1 & 1 & 0
\end{matrix}
\right) \]
(up to row and column swaps).
\end{lemma}

\section{Number of monomial Cremona transformations}
\label{sec:num_of_transformations}

\begin{thm}
\label{thm:num_of_transformations}
For $\delta \geq 3$, there are exactly $\phi(\delta)$ exponent matrices, as 
described in \autoref{thm:mainThm}, where $\phi$ is Euler's totient function.

\begin{proof}
There is exactly one exponent matrix of the form
\[ \left( \begin{matrix}
			0 & 0 & \delta\\
			1 & \delta-1 & 0\\
			0 & 1 & \delta -1
\end{matrix} \right), \]
so we must show that there are $\phi(\delta)-1$ exponent matrices of the form 
described by ({\bf III\/}) in \autoref{thm:mainThm}.  These are
\[ \left(
		\begin{matrix}
			0 & 0 & \delta \\
			a & b & 0 \\
			c & d & e
		\end{matrix}
		\right) \]
		satisfying
		\begin{enumerate}[\rm(i)]
			\item $a+b = \delta$,\; $c+d+e=\delta$,
			\item $ad-bc=1$.
		\end{enumerate}
Substituting $d=\delta-c-e$ and $b=\delta-a$ into (iii), we get
\begin{equation} \label{eqn:thm3_vol}
	\delta a - \delta c - ae = 1
\end{equation}
which we rearrange to get 
\[\delta(a-c) = ae + 1.\]
This implies $\delta \;\big|\; (ae + 1)$ which we can write as 
$ae \equiv \delta-1 \mod{\delta}$.  Then, for each $a \in \bbZ/\delta\bbZ$ we 
can set $e = a^{-1} (\delta-1)$.  However, we cannot allow $e = \delta-1$ 
because this would require $c+d=1$, which is impossible by $(ii)$.  

If $a,e$ are determined, then $\delta r = ae +1$ for some positive 
integer $r$, and we can let $c = r-a$.

Thus, there is a unique solution $\{a,c,e\}$ to \eqref{eqn:thm3_vol} for each 
value of $a \in \bbZ/\delta\bbZ \backslash \{1\}$.  Since this set has order 
$\phi(\delta)$-1, the proof is complete.
\end{proof}
\end{thm}

By \autoref{thm:mainThm}, there are exactly 2 exponent matrices having 
$\delta =2$.  These are
\[
\left( 
\begin{matrix}
0 & 1 & 1 \\
1 & 0 & 1 \\
1 & 1 & 0
\end{matrix}
\right), \quad
\left( 
\begin{matrix}
0 & 0 & 2 \\
1 & 1 & 0 \\
0 & 1 & 1
\end{matrix}
\right)
\]
and as always, this statement is up to row/column swapping.  

We could consider these 
maps up to \emph{permutation similarity}, where two exponent matrices $M,M'$ 
yield the same map if and only if there exists a permutation matrix $Z$ such 
that $M' = Z^\textrm{T} M Z$.  In other words, $M,M'$ are equivalent if and 
only if there exists $\sigma \in S_3$, the symmetric group on $\{1,2,3\}$, such 
that $M = \sigma \star M'$ where $S_3$ acts on $M_3(\mathbb{Z})$ by permuting 
the rows and columns of $M'$.  For instance, if $\sigma = (1\;2)$, then $\sigma 
\star M'$ is the matrix obtained from $M'$ by swapping the first and second 
rows and then swapping the first and second columns. 

Since $\#(S_3) = 6$ and $S_3$ acts freely on the subset of $M_3(\mathbb{Z})$ 
defined by \autoref{thm:mainThm}, under this new equivalence relation we 
have exactly 6 non-similar matrices for each matrix described by the theorem.
So in the new context, we say that there are exactly 12 
monomial Cremona transformations on $\bbP^2$ having $\delta = 2$.

Under this equivalence relation we have the 
\begin{corollary}
For $\delta \geq 3$ there are exactly $6 \cdot \phi(\delta)$ monomial Cremona 
maps on $\bbP^2$ up to permutation similarity.
\qed
\end{corollary}

To motivate this perspective, consider the problem of determining $\gamma_1 = \delta$ for $\varphi^n = \varphi \circ \varphi \circ \cdots \circ \varphi$ for some monomial Cremona transformation $\varphi$.  Let $\varphi_1, \varphi_2$ be the two maps defined by 
\[ M_1 = \left(\begin{matrix}
0 & 0 & 5 \\
4 & 1 & 0 \\
3 & 1 & 1
\end{matrix}
\right), \quad
M_2 = \left( 
\begin{matrix}
3 & 1 & 1 \\
0 & 0 & 5 \\
4 & 1 & 0
\end{matrix} \right), \]
that is, $M_2$ is given by acting with $(1 2 3)$ on the rows of $M_1$,
Then, letting $d_n$ denote $\delta = \gamma_1$ corresponding to $\varphi^n$, we 
have 
\[\begin{array}{c|c c c c c}
	& d_1 & d_2 & d_3 & d_4 & d_5 \\ \hline
	\varphi_1 & 5 & 15 & 40 & 105 & 275 \\
	\varphi_2 & 5 & 9 & 13 & 17 & 21
\end{array}\]
so these two maps should not be considered `the same'.

\bibliography{refs}{}
\bibliographystyle{amsalphaeprint}


\end{document}